\documentclass[12pt]{amsproc}

\usepackage[all]{xy}
\usepackage{latexsym,amsthm,amsmath,amssymb,enumerate}

\newtheorem{thm}{Theorem}[section]
\newtheorem{prop}[thm]{Proposition}

\theoremstyle{definition}
\newtheorem{definition}[thm]{Definition}
\newtheorem{exa}[thm]{Example}
\newtheorem{remark}[thm]{Remark}

\begin{document}

\title[Componentwise linearity of ideals]{Componentwise linearity of ideals\\ arising from graphs}

\author{Veronica Crispin}
\thanks{The first author was sponsored by The Royal Swedish Academy of Sciences}
\address{Department of Mathematics\\
University of Oregon}
\email{vcrispin@uoregon.edu}

\author{Eric Emtander}
\address{Department of Mathematics\\
Stockholm University}
\email{erice@math.su.se}

\keywords{Linear resolutions, componentwise linear, chordal graphs, monomial ideals, polymatroidal ideals}
\subjclass[2000]{Primary 13C05, 05E40; Secondary 13A15}

\begin{abstract}
Let $G$ be a simple undirected graph on $n$ vertices.  Francisco and Van Tuyl 
have shown that if $G$ is chordal, then $\bigcap_{\{ x_i,x_j\}\in E_G} 
\langle x_i,x_j\rangle$ is componentwise linear. A natural
question that arises is for which $t_{ij}>1$ the ideal 
$\bigcap_{\{x_i,x_j\}\in E_G}\langle x_i, x_j\rangle^{t_{ij}}$ is
componentwise linear, if $G$ is chordal. In this report we show that
$\bigcap_{\{ x_i,x_j\}\in E_G} \langle x_i, x_j\rangle^{t}$ is componentwise linear
for all $n\geq 3$ and positive $t$, if $G$ is a complete graph. We give also an example
where $G$ is chordal, but the intersection ideal is not componentwise
linear for any $t>1$.
\end{abstract}

\maketitle

\section{Introduction}
The previous version of this article is published as \cite{CQE}.

\bigskip
Let $G$ be a simple graph on $n$ vertices, $E_G$ the edge set of $G$
and $V_G$ the vertex set of $G$. Let $R=k[x_1, \ldots, x_n]$ be the polynomial ring over a field $k$. The {\em edge ideal} of $G$ is
the quadratic squarefree monomial ideal ${\mathcal
  I}(G)=\langle\{ x_ix_j\}\;\vert\; \{ x_i, x_j\}\in E_G\rangle\subset R$. Then we
define the {\em squarefree Alexander dual} of ${\mathcal I}(G)$ as
${\mathcal I}(G)^\vee=\cap_{\{ x_i,x_j\}\in E_G} \langle x_i, x_j\rangle$. Calling 
${\mathcal I}(G)^\vee$ the squarefree Alexander dual of ${\mathcal I}(G)$ is natural since ${\mathcal I}(G)^\vee$ 
is the Stanley--Reisner ideal of the simplicial complex $\Delta^\vee$, that is, the 
Alexander dual simplicial complex of $\Delta$. Here $\Delta$ is the simplicial 
complex, the Stanley-Reisner ideal of which is ${\mathcal I}(G)$. 

In \cite{HH} Herzog and Hibi give the following definition. Given a graded ideal 
$I\subset R$, we denote by $I_{\langle d\rangle}$ the ideal generated by the elements 
of  degree $d$ that belong to $I$. Then we say that a (graded) ideal  $I\subset R$ is 
{\em componentwise linear} if $I_{\langle d\rangle}$ has a linear resolution for all $d$.

If the graph $G$ is chordal, that is, every cycle of length $m\geq 4$ in $G$ has a chord, 
then it is proved by Francisco and Van Tuyl \cite{FvT1} that ${\mathcal I}(G)^V$ is 
componentwise linear. (The authors then use the result to show that all chordal graphs are 
sequentially Cohen-Macaulay.)

In this report we examine componentwise linearity of ideals arising
from complete graphs and of the form
$\bigcap_{\{ x_i,x_j\}\in E_G}\langle x_i, x_j\rangle^{t}$.

\section{Intersections for complete graphs}\label{iocg}
Let $K_n$ be a complete graph on $n$ vertices, that is, $\{
x_i,x_j\}\in E_{K_n}$ for all $1\leq i\ne j\leq n$. We write
$K_n^{(t)}=\bigcap_{\{ x_i,x_j\}\in E_{K_n}}\langle x_i,
x_j\rangle^{t}$. We will show that the ideal $K_n^{(t)}$ is
componentwise linear for all $n\geq 3$ and $t\geq 1$. Recall that a {\em
  vertex cover} of a graph $G$ is a subset $A\subset V_G$ such that
every edge of $G$ is incident to at least one vertex of $A$. One can
show that $\mathcal I(G)^V=\langle x_{i_1}\cdots x_{i_k}\;\vert\;
\{x_{i_1},\ldots ,x_{i_k}\}\; {\rm a\; vertex\; cover\;of}\; G\rangle$. A
{\em t-vertex cover} (or a {\em vertex cover of order t) of G} is a
vector ${\rm\bf a} =(a_1,\ldots, a_n)$ with $a_i\in\mathbb N$ such that 
$a_i+a_j\geq t$ for all $\{ x_i,x_j\}\in E_G$.

In the proof of our main result Theorem 2.3, we use the following definition and proposition.

\begin{definition}\label{lq} A monomial ideal $I$ is said to
have {\em linear quotients}, if for some degree ordering of the minimal
generators $f_1,\ldots ,f_r$ and all $k>1$, the colon ideals $\langle f_1,\ldots ,f_{k-1}\rangle :f_k$ are generated by a subset of $\{x_1,\ldots ,x_n \}$.
\end{definition}

\begin{prop}[Proposition 2.6 in \cite{FvT2} and Lemma 4.1 in \cite{CH}]\label{lqcl}
If $I$ is a homogeneous ideal with linear quotients, then $I$ is
componentwise linear.
\end{prop}

\begin{thm}
The ideal $K_n^{(t)}$ is componentwise linear for all $n\geq 3$ and $t\geq 1$.

\begin{proof} For calculating an explicit generating system of $K_n^{(t)}$ 
we will use $t$-vertex covers. Pick any monomial $m$ in the generating set of $K_n^{(t)}$ and, 
for some $k$ and $l$, consider the greatest exponents $t_k$ and $t_l$ such that $x_k^{t_k}x_l^{t_l}$ is 
a factor in $m$. As $m$ is contained in $\langle x_k,x_l\rangle^{t}$ we must 
have $t_k+t_l\geq t$. Hence, $K_n^{(t)}$ is generated by the monomials of the form 
$\rm\bf x^a$, where $\rm\bf a$ is an $t$-cover of $K_n$. That is, the sum of the two 
lowest exponents in every (monomial) generator of $K_n^{(t)}$ is at least $t$.

First we assume that $t=2m+1$ is odd. Using the degree lexicographic ordering 
$x_1\prec x_2\prec\cdots \prec x_n$ on the the minimal generators we get
\begin{displaymath}
\begin{array}{lrcr}
  K_n^{(t)}=K_n^{(2m+1)}= & \big\langle x_1^{m}\prod_{i\ne 1} x_i^{m+1}, &
  \ldots & ,x_n^{m}\prod_{i\ne n} x_i^{m+1}, \\
&&&\\
  & x_1^{m-1}\prod_{i\ne 1} x_i^{m+2}, &\ldots & ,x_n^{m-1}\prod_{i\ne
   1} x_i^{m+2},\\
& & \vdots & \\
 & \prod_{i\ne 1} x_i^{2m+1}, &\ldots & ,\prod_{i\ne
   n} x_i^{2m+1}\big\rangle. 
\end{array} 
\end{displaymath}
This ordering of the minimal generators satisfies the condition in
Definition~\ref{lq}. Hence, $K_n^{(t)}$ has linear quotients and is
componentwise linear by Proposition~\ref{lqcl}. 

If $t=2m$ is even, then the degree lexicographic ordering yields the sequence
\begin{displaymath}
\begin{array}{lrcr}
  K_n^{(t)}=K_n^{(2m)}= & \big\langle \prod_{i=1}^{2m} x_i^m,\quad  x_1^{m-1}\prod_{i\ne 1} x_i^{m+1}, &
  \ldots & ,x_n^{m-1}\prod_{i\ne n} x_i^{m+1}, \\
&&&\\
  & x_1^{m-2}\prod_{i\ne 1} x_i^{m+2}, &\ldots & ,x_n^{m-2}\prod_{i\ne
   1} x_i^{m+2},\\
& & \vdots & \\
 & \prod_{i\ne 1} x_i^{2m}, &\ldots & ,\prod_{i\ne
   n} x_i^{2m}\big\rangle, 
\end{array} 
\end{displaymath}
which also satisfies the condition in Definition~\ref{lq}, and the same result follows. 
\end{proof}
\end{thm}

\begin{exa}
\[
K_{12}^{(5)}=\big\langle \{ x_j^2\prod_{i\ne j}
x_i^3\}_{1\leq j\leq 12},\; \{ x_j\prod_{i\ne j}
x_i^4\}_{1\leq j\leq 12},\;
\{\prod_{i\ne j} x_i^5\}_{1\leq j\leq 12} \big\rangle
\]
and
\[
K_5^{(6)}=\big\langle \prod_{i=1}^5 x_i^3,\;\{ x_j^2\prod_{i\ne j}
x_i^4\}_{1\leq j\leq 5},\;\{ x_j\prod_{i\ne j}
x_i^5\}_{1\leq j\leq 5},\;\{\prod_{i\ne j} x_i^6\}_{1\leq j\leq 5}
\big\rangle.
\]
\end{exa}

\begin{remark}
A monomial ideal is called $polymatroidal$ if it is generated in one degree and its minimal generators satisfy a certain "exchange condition". In \cite{HT} Herzog and Takayama show that polymatroidal ideals have linear resolutions. Later Francisco and van Tuyl \cite{FvT2} proved that  some families of ideals $I$ are componentwise linear showing in their Theorem~3.1 that $I_{\langle d\rangle}$ are polymatroidal for all $d$.

The ideals $K_n^{(t)}$ are also polymatroidal, but the proof using the same techniques as in the proof of Theorem~3.1 in \cite{FvT2}  is rather tedious and takes a few pages.

\end{remark}

\section{A counterexample }

{\em There exists a chordal graph $G$ such that 
$\bigcap_{\{ x_i,x_j\}\in E_G}\langle x_i, x_j\rangle^t$ is not
componentwise linear for any $t>1$.}

\begin{proof}
Let $G$ be the chordal graph
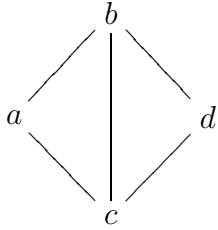
\[
\xymatrix{&b\ar@{-}[dl]\ar@{-}[dd]\ar@{-}[dr]&\\
a\ar@{-}[dr]&&d\ar@{-}[dl]\\
&c&}
\]
and denote the intersection $\langle a,b\rangle^t\cap\langle a,c\rangle^t\cap\langle b,c\rangle^t\cap\langle b,d\rangle^t  $ by $I_4^{(t)}$. We have
\[
I_4^{(1)}=\bigcap_{\{ i,j \}\in E_G}\langle i,j\rangle=\langle bc \rangle + \langle abd,acd\rangle
\]
and
\[
I_4^{(2)}=\bigcap_{\{ i,j \}\in E_G}\langle i,j\rangle =\langle b^2c^2,abcd\rangle + \langle a^2b^2d^2,a^2c^2d^2\rangle.
\]

Arguing in the same way as for $K_n^{(t)}$ we see that the minimal generating set consists of generators of exactly degree $2t$ and generators of higher degrees:

\begin{itemize}
\item If $t_a\leq\lfloor\frac{t}{2} \rfloor$ then $t_b=t-t_a=t_c$ (the sum $t_b+t_c\geq t$ automatically) and $t_d=t-t_b=t-t_c=t_a$. We get the set of minimal generators of degree $2t$:
\[
\big\{ a^{i}(bc)^{t-i}d^{i} \big\}_{0\leq i\leq\lfloor\frac{t}{2} \rfloor }.
\]

\item If $t_a> \lfloor\frac{t}{2} \rfloor$, then either $t_b=t-t_a$ and $t_c=t-t_b=t_a$, or $t_c=t-t_a$ and $t_b=t_a$. Further $t_d=t_a$. The set of minimal generators we get in this way is equal to
\[
\big\{ (acd)^{i}b^{t-i} \big\}_{\lfloor\frac{t}{2} \rfloor < i\leq t} \cup \big\{ (acd)^{i}b^{t-i} \big\}_{\lfloor\frac{t}{2} \rfloor  < i\leq t }.
\]
The generators in this set are of degree at least $(2t+1)$ for odd $t$ and of degree at least $(2t+2)$ for even $t$.
\end{itemize}

\bigskip

Now consider the minimal free resolution $\mathcal{F}.$ of $(I_4^{(t)})_{\langle 2t\rangle}$. Since $\mathcal{F}.$ 
is contained in any free resolution $\mathcal{G}.$ of $(I_4^{(t)})_{\langle 2t\rangle}$ we have that if $F_1$ 
(the component of $\mathcal{F}$. in homological degree 1) has a non-zero component in a certain degree, then so 
does $G_1$. Let $\mathcal{G}.$ be the Taylor resolution of $(I_4^{(t)})_{\langle 2t\rangle}$. The degrees in which 
$G_1$ has nonzero components come from least common mutliples of pairs of minimal generators of 
$(I_4^{(t)})_{\langle 2t\rangle}$. By considering the above description of the minimal generators in degree $2t$, 
one sees that $G_1$ has non-zero components only in degrees strictly larger than $2t+1$. Thus $\mathcal{F}.$ cannot 
be a linear resolution
and, hence, $I_4^{(t)}$ is not componentwise linear.
\end{proof}



\section*{Acknowledgements}
First of all we would like to thank the Universit$\rm\grave{a}$ di Catania and the organizers
 of the PRAGMATIC summer school 2008, especially Alfio Ragusa and Giuseppe Zappal$\rm\grave{a}$. We are deeply
 greatful to J\"urgen Herzog and Volkmar Welker for their excellent
 lectures, interesting problems and thorough guidance.

\bibliographystyle{amsalpha}

\end{document}